\documentstyle{amsppt}
\voffset-12mm
\magnification1200
\pagewidth{130mm}
\pageheight{204mm}
\hfuzz=2.5pt\rightskip=0pt plus1pt
\binoppenalty=10000\relpenalty=10000\relax
\TagsOnRight
\loadbold
\nologo
\addto\tenpoint{\normalbaselineskip=1.3\normalbaselineskip\normalbaselines}
\addto\eightpoint{\normalbaselineskip=1.2\normalbaselineskip\normalbaselines}
%==================================================
\let\le\leqslant
\let\ge\geqslant
\let\<\langle
\let\>\rangle
\let\epsilon\varepsilon
\redefine\d{\roman d}
\redefine\Re{\operatorname{Re}}
\define\ba{\boldkey a}
\define\bb{\boldkey b}
\define\bh{\boldkey h}
\define\fc{\frak c}
\define\fG{\frak G}
%==================================================
\topmatter
\title
Multiple-integral representations \\
of very-well-poised hypergeometric series
\footnotemark
\endtitle
\author
Wadim Zudilin \rm(Moscow)
\endauthor
\date
\hbox to100mm{\vbox{\hsize=100mm%
\centerline{E-print \tt math.CA/0206177}
\smallskip
\centerline{17 March 2002}
}}
\enddate
\address
\hbox to70mm{\vbox{\hsize=70mm%
\leftline{Moscow Lomonosov State University}
\leftline{Department of Mechanics and Mathematics}
\leftline{Vorobiovy Gory, GSP-2}
\leftline{119992 Moscow, RUSSIA}
\leftline{{\it URL\/}: \tt http://wain.mi.ras.ru/index.html}
}}
\endaddress
\email
{\tt wadim\@ips.ras.ru}
\endemail
\endtopmatter
\footnotetext{An extract from my contribution
``Well-poised hypergeometric service
for diophantine problems of zeta values''
to \it Actes des 12\`emes rencontres arithm\'etiques de Caen.}
\rightheadtext{Hypergeometric series and multiple integrals}
\leftheadtext{W.~Zudilin}
%==================================================
\document

The aim of this note is to connect two objects:
very-well-poised hypergeometric series
$$
\align
F_k(\bh)
&=F_k(h_0;h_1,\dots,h_k)
:=\frac{\Gamma(1+h_0)\,\prod_{j=1}^k\Gamma(h_j)}
{\prod_{j=1}^k\Gamma(1+h_0-h_j)}
\\ &\qquad\times
{}_{k+2}\!F_{k+1}\biggl(\matrix\format&\,\c\\
h_0, & 1+\tfrac12h_0, & h_1, & \dots, & h_k \\
& \tfrac12h_0, & 1+h_0-h_1, & \dots, & 1+h_0-h_k
\endmatrix\biggm|(-1)^{k+1}\biggr)
\\ &\phantom:
=\sum_{\mu=0}^\infty(h_0+2\mu)
\frac{\prod_{j=0}^k\Gamma(h_j+\mu)}{\prod_{j=0}^k\Gamma(1+h_0-h_j+\mu)}
(-1)^{(k+1)\mu},
\tag1
\endalign
$$
and multiple integrals
$$
\align
J_k(\ba,\bb)
&\phantom:=J_k\biggl(\matrix\format&\,\c\\
a_0, & a_1, & \dots, & a_k \\
& b_1, & \dots, & b_k
\endmatrix\biggr)
\\
&:=\idotsint\limits_{[0,1]^k}
\frac{\prod_{j=1}^kx_j^{a_j-1}(1-x_j)^{b_j-a_j-1}}
{Q_k(x_1,x_2,\dots,x_k)^{a_0}}
\,\d x_1\,\d x_2\dotsb\d x_k
\tag2
\endalign
$$
with $Q_0:=1$ and
$$
\alignat1
Q_k
&=Q_k(x_1,x_2,\dots,x_k)
:=1-(1-(\dotsb(1-(1-x_k)x_{k-1})\dotsb)x_2)x_1
\\
&=1-x_1Q_{k-1}(x_2,\dots,x_k)
=Q_{k-1}(x_1,\dots,x_{k-1})+(-1)^kx_1x_2\dotsb x_k
\tag3
\endalignat
$$
for $k\ge1$.

\proclaim{Theorem}
For each $k\ge1$, there holds the identity
$$
\alignat1
&
\frac{\prod_{j=1}^{k+1}\Gamma(1+h_0-h_j-h_{j+1})}
{\Gamma(h_1)\,\Gamma(h_{k+2})}
\cdot F_{k+2}(h_0;h_1,\dots,h_{k+2})
\\ &\qquad
=J_k\biggl(\matrix\format&\,\c\\
h_1, & h_2, & h_3, & \dots, & h_{k+1} \\
& 1+h_0-h_3, & 1+h_0-h_4, & \dots, & 1+h_0-h_{k+2}
\endmatrix\biggr),
\tag4
\endalignat
$$
provided that
$$
\gather
1+\Re h_0>\frac2{k+1}\cdot\sum_{j=1}^{k+2}\Re h_j,
\tag5
\\
\Re(1+h_0-h_{j+1})>\Re h_j>0
\qquad\text{for}\quad j=2,\dots,k+1,
\tag6
\\ \vspace{1.5pt}
h_1,h_{k+2}\ne0,-1,-2,\dots\,.
\tag7
\endgather
$$
\endproclaim

\remark{Remark}
Condition~\thetag{5} is required for the absolute convergence
of the series~\thetag{1} in the unit circle (and, in particular,
at the point~$(-1)^{k+1}$), while condition~\thetag{6} ensures
the convergence of the corresponding multiple integral~\thetag{2}.
The restriction~\thetag{7} can be removed by the theory of
analytic continuation if we write
$\Gamma(h_j+\mu)/\Gamma(h_j)$ for $j=1,k+2$ as Pochhammer's
symbol $(h_j)_\mu$ when summing in~\thetag{1}.
\endremark

\smallskip
In the case of {\it integral\/} parameters~$\bh$,
the quantities~\thetag{1} are known to be $\Bbb Q$-linear
forms in even/odd zeta values depending on evenness/oddness
of $k\ge4$. Therefore, if positive integral parameters
$\ba$ and $\bb$ satisfy the additional condition
$$
b_1+a_2=b_2+a_3=\dots=b_{k-1}+a_k,
\tag8
$$
then the quantities~\thetag{2} are $\Bbb Q$-linear
forms in even/odd zeta values.
Specialization $a_j=n+1$ and $b_j=2n+2$ gives one the coincidence
of multiple integrals and well-poised hypergeometric series
conjectured by us in~\cite{Zu2, Section~9}.
Denoting the corresponding integrals by~$J_{k,n}$
and using our arithmetic results~\cite{Zu1, Lemmas~4.2--4.4}
we then conclude that
$$
D_n^{k+1}\Phi_n^{-1}J_{k,n}
\in\Bbb Z\zeta(k)+\Bbb Z\zeta(k-2)+\dots+\Bbb Z\zeta(3)+\Bbb Z
\quad\text{for $k$ odd},
$$
where
$$
\Phi_n:=\prod\Sb p<n\\\{n/p\}\in[\frac23,1)\endSb p,
\qquad
\lim_{n\to\infty}\frac{\log\Phi_n}n
=\psi(1)-\psi\Bigl(\frac23\Bigr)-\frac12
=0.24101875\dots
$$
($\{\,\cdot\,\}$ denotes the fractional part of a number),
that is closed enough to Vasilyev's conjectural
inclusions~\cite{V}.
The choice $a_j=rn+1$ and $b_j=(r+1)n+2$ in~\thetag{2}
(or, equivalently, $h_0=(2r+1)n+2$ and $h_j=rn+1$
for $j=1,\dots,k+2$ in~\thetag{1})
with the integer $r\ge1$
depending on a given odd integer~$k$ presents
almost the same linear forms in odd zeta values
as considered by T.~Rivoal in~\cite{R} for proving
his remarkable result on infiniteness of irrational
numbers in the set $\zeta(3),\zeta(5),\zeta(7),\dots$\,.

In addition, we have to mention, under hypothesis~\thetag{8},
the obvious stability of the quantity
$$
\align
\frac{F_{k+2}(h_0;h_1,\dots,h_{k+2})}
{\prod_{j=1}^{k+2}\Gamma(h_j)}
&=\frac{J_k(\ba,\bb)}
{\prod_{j=2}^{k+1}\Gamma(h_j)\cdot\prod_{j=1}^{k+1}\Gamma(1+h_0-h_j-h_{j+1})}
\\
&=\frac{J_k(\ba,\bb)}
{\prod_{j=1}^k\Gamma(a_j)\cdot\Gamma(b_1+a_2-a_0-a_1)
\cdot\prod_{j=1}^k\Gamma(b_j-a_j)}
\endalign
$$
under the action of the ($\bh$-trivial) group~$\fG$
of order~$(k+2)!$ containing all permutations
of the parameters $h_1,\dots,h_{k+2}$.
This fact can be applied for number-theoretical applications.
In the cases $k=2$ and $k=3$ the change of variables
$(x_{k-1},x_k)\mapsto(1-x_k,1-x_{k-1})$ in~\thetag{2}
produces an additional transformation~$\fc$ of both~\thetag{2}
and~\thetag{1}; for $k\ge4$ this transformation is not yet
available since condition~\thetag{8} is broken.
The groups $\<\fG,\fc\>$ of orders~$120$ and~$1920$
for $k=2$ and $k=3$ respectively are known~\cite{RV1},~\cite{RV2};
G.~Rhin and C.~Viola make a use of these groups to discover
nice estimates for the irrationality measures of~$\zeta(2)$ and~$\zeta(3)$.
Finally, we want to note that the group~$\fG$ can be easily
interpretated as the permutation group of the parameters
$$
e_{0l}=h_l-1, \quad 1\le l\le k+2,
\qquad e_{jl}=h_0-h_j-h_l, \quad 1\le j<l\le k+2
$$
(see~\cite{Zu2, Section~9} for details).

\proclaim{Lemma 1}
Theorem is true if $k=1$.
\endproclaim

\demo{Proof}
Thanks to a limiting case of Dougall's theorem,
$$
F_3(h_0;h_1,h_2,h_3)
=\frac{\Gamma(h_1)\,\Gamma(h_2)\,\Gamma(h_3)\,\Gamma(1+h_0-h_1-h_2-h_3)}
{\Gamma(1+h_0-h_1-h_2)\,\Gamma(1+h_0-h_1-h_3)\,\Gamma(1+h_0-h_2-h_3)}
\tag9
$$
(see, e.g., \cite{B, Section~4.4, formula~(1)}), provided that
$1+\Re h_0>\Re(h_1+h_2+h_3)$ and $h_j$~is not a non-positive
integer for $j=1,2,3$. On the other hand, the integral on the right
of~\thetag{4} has Euler type, that is
$$
\align
J_1\biggl(\matrix\format&\,\c\\
h_1, & h_2 \\ & 1+h_0-h_3
\endmatrix\biggr)
&=\int_0^1\frac{x^{h_2-1}(1-x)^{h_0-h_2-h_3}}{(1-x)^{h_1}}\,\d x
\\
&=\frac{\Gamma(h_2)\,\Gamma(1+h_0-h_1-h_2-h_3)}{\Gamma(1+h_0-h_1-h_3)},
\endalign
$$
provided that $1+\Re h_0>\Re(h_1+h_2+h_3)$ and $\Re h_2>0$.
Therefore, multiplying equality~\thetag{9} by the required
product of gamma-functions we deduce identity~\thetag{4}
if $k=1$.
\enddemo

\remark{Remark}
If we arrange about $J_0(a_0)$ to be~$1$,
the claim of Theorem remains valid if $k=0$
thanks to another consequence of Dougall's theorem
\cite{B, Section~4.4, formula~(3)}.
\endremark

\proclaim{Lemma 2 \cite{N, Section~3.2}}
Let $a_0,a,b\in\Bbb C$ and $t_0\in\Bbb R$ be numbers
satisfying the conditions
$$
\Re a_0>t_0>0, \quad \Re a>t_0>0,
\quad\text{and}\quad \Re b>\Re a_0+\Re a.
$$
Then for any non-zero $z\in\Bbb C\setminus(1,+\infty)$
the following identity holds:
$$
\align
&
\int_0^1\frac{x^{a-1}(1-x)^{b-a-1}}{(1-zx)^{a_0}}\,\d x
\\ &\qquad
=\frac{\Gamma(b-a)}{\Gamma(a_0)}
\cdot\frac1{2\pi i}\int_{-t_0-i\infty}^{-t_0+i\infty}
\frac{\Gamma(a_0+t)\,\Gamma(a+t)\,\Gamma(-t)}{\Gamma(b+t)}
\,(-z)^t\,\d t,
\tag10
\endalign
$$
where $(-z)^t=|z|^te^{it\arg(-z)}$, $-\pi<\arg(-z)<\pi$
for $z\in\Bbb C\setminus[0,+\infty)$ and $\arg(-z)=\pm\pi$
for $z\in(0,1]$. The integral on the right-hand side of~\thetag{10}
converges absolutely. In addition, if $|z|\le1$, both integrals
in~\thetag{10} can be identified with the absolutely convergent
Gauss hypergeometric series
$$
\frac{\Gamma(a)\,\Gamma(b-a)}{\Gamma(b)}
\cdot{}_2\!F_1\biggl(\matrix\format&\,\c\\
a_0, & a \\ & b \endmatrix\biggm|z\biggr)
=\frac{\Gamma(b-a)}{\Gamma(a_0)}\sum_{\nu=0}^\infty
\frac{\Gamma(a_0+\nu)\,\Gamma(a+\nu)}{\nu!\,\Gamma(b+\nu)}z^\nu.
$$
\endproclaim

Set $\epsilon_k=0$ for $k$~even
and $\epsilon_k=1$ or~$-1$ for $k$~odd.

\proclaim{Lemma 3}
For each integer $k\ge2$, there holds the relation
$$
\align
J_k\biggl(\matrix\format&\,\c\\
a_0, & a_1, & \dots, & a_{k-1}, & a_k \\
& b_1, & \dots, & b_{k-1}, & b_k
\endmatrix\biggr)
&=\frac{\Gamma(b_k-a_k)}{\Gamma(a_0)}
\cdot\frac1{2\pi i}\int_{-t_0-i\infty}^{-t_0+i\infty}
\frac{\Gamma(a_0+t)\,\Gamma(a_k+t)\,\Gamma(-t)}{\Gamma(b_k+t)}
\\ &\qquad\times
e^{\epsilon_k\pi it}\cdot J_{k-1}\biggl(\matrix\format&\,\c\\
a_0+t, & a_1+t, & \dots, & a_{k-1}+t \\
& b_1+t, & \dots, & b_{k-1}+t
\endmatrix\biggr)\,\d t,
\endalign
$$
provided that $\Re a_0>t_0>0$, $\Re a_k>t_0>0$,
$\Re b_k>\Re a_0+\Re a_k$, and the integral on the left
converges.
\endproclaim

\demo{Proof}
We start with mentioning that the first recursion in~\thetag{3}
and inductive arguments yield the inequality
$$
0<Q_k(x_1,x_2,\dots,x_k)<1 \qquad\text{for}\quad
(x_1,x_2,\dots,x_k)\in(0,1)^k \quad\text{and}\quad k\ge1.
\tag11
$$
By the second recursion in~\thetag{3},
\ $Q_k=Q_{k-1}\cdot(1-zx_k)$ for $k\ge2$, where
$$
z=\frac{(-1)^{k+1}x_1\dotsb x_{k-1}}{Q_{k-1}(x_1,\dots,x_{k-1})}.
$$
For each $(x_1,\dots,x_{k-1})\in(0,1)^{k-1}$,
the number~$z$ is real with the property $z<0$ for $k$~even
and $0<z<1$ for $k$~odd, since in the last case we have
$$
z=\frac{x_1\dotsb x_{k-1}}{Q_{k-1}(x_1,\dots,x_{k-2},x_{k-1})}
=\frac{x_1\dotsb x_{k-1}}{Q_{k-2}(x_1,\dots,x_{k-2})+x_1\dotsb x_{k-1}}
<1
$$
by~\thetag{11}. Therefore, splitting the integral~\thetag{2}
over $[0,1]^k=[0,1]^{k-1}\times[0,1]$ and applying Lemma~2 to the integral
$$
\int_0^1\frac{x_k^{a_k-1}(1-x_k)^{b_k-a_k-1}}{(1-zx_k)^{a_0}}\,\d x_k
$$
we arrive at the desired relation.
\enddemo

\demo{Proof of Theorem}
The case $k=1$ is considered in Lemma~1.
Therefore we will assume that $k\ge2$, identity~\thetag{4}
holds for~$k-1$, and, in addition, that
$$
1+\Re h_0>\frac 2k\cdot\sum_{j=1}^{k+1}\Re h_j,
\qquad \Re h_{k+2}<1.
\tag12
$$
The restrictions~\thetag{12} can be easily removed from the
final result by the theory of analytic continuation.

By the inductive hypothesis, for $t\in\Bbb C$ with $\Re t<0$,
we deduce that
$$
\align
&
J_{k-1}\biggl(\matrix\format&\,\c\\
h_1+t, & h_2+t, & h_3+t, & \dots, & h_k+t \\
& 1+h_0-h_3+t, & 1+h_0-h_4+t, & \dots, & 1+h_0-h_{k+1}+t
\endmatrix\biggr)
\\ &\qquad
=\frac{\prod_{j=1}^k\Gamma(1+h_0-h_j-h_{j+1})}
{\Gamma(h_1+t)\,\Gamma(h_{k+1}+t)}
\cdot F_{k+1}(h_0+2t;h_1+t,\dots,h_{k+1}+t)
\\ &\qquad
=\frac{\prod_{j=1}^k\Gamma(1+h_0-h_j-h_{j+1})}
{\Gamma(h_1+t)\,\Gamma(h_{k+1}+t)}
\cdot\frac1{2\pi i}\int_{-s_0-i\infty}^{-s_0+i\infty}
(h_0+2t+2s)
\\ &\qquad\quad\times
\frac{\Gamma(h_0+2t+s)\,\prod_{j=1}^{k+1}\Gamma(h_j+t+s)\,\Gamma(-s)}
{\prod_{j=1}^{k+1}\Gamma(1+h_0-h_j+t+s)}
\,e^{\epsilon_{k-1}\pi is}\,\d s,
\tag13
\endalign
$$
where the real number $s_0>0$ satisfies the conditions
$$
\Re(h_0+2t)>s_0, \quad \Re(1+\tfrac12h_0+t)>s_0,
\qquad \Re(h_j+t)>s_0 \quad\text{for}\; j=1,\dots,k+1,
$$
and the absolute convergence of the last Barnes-type
integral follows from~\cite{N, Lem\-ma~3}.
Shifting the variable $t+s\mapsto s$ in~\thetag{13}
(with a help of the equality
$e^{\epsilon_k\pi it}\cdot\allowmathbreak e^{\epsilon_{k-1}\pi is}
=e^{\epsilon_{k-1}\pi i(t+s)}\cdot e^{\epsilon_1\pi it}$),
applying Lemma~3, and interchanging double integration
(thanks to the absolute convergence of the integrals)
we conclude that
$$
\align
&
J_k\biggl(\matrix\format&\,\c\\
h_1, & h_2, & h_3, & \dots, & h_k, & h_{k+1} \\
& 1+h_0-h_3, & 1+h_0-h_4, & \dots, & 1+h_0-h_{k+1}, & 1+h_0-h_{k+2}
\endmatrix\biggr)
\\ &\qquad
=\frac{\prod_{j=1}^{k+1}\Gamma(1+h_0-h_j-h_{j+1})}{\Gamma(h_1)}
\\ &\qquad\quad\times
\frac1{2\pi i}\int_{-s_1-i\infty}^{-s_1+i\infty}
(h_0+2s)\frac{\prod_{j=1}^{k+1}\Gamma(h_j+s)}
{\prod_{j=1}^{k+1}\Gamma(1+h_0-h_j+s)}e^{\epsilon_{k-1}\pi is}
\\ &\qquad\quad\times
\frac1{2\pi i}\int_{-t_0-i\infty}^{-t_0+i\infty}
\frac{\Gamma(-s+t)\,\Gamma(h_0+s+t)\,\Gamma(-t)}
{\Gamma(1+h_0-h_{k+2}+t)}
e^{\epsilon_1\pi it}\,\d t\,\d s,
\tag14
\endalign
$$
where $s_1=s_0+t_0$.
Since $\Re h_{k+2}<1$ and $h_{k+2}\ne0,-1,-2,\dots$,
the last Barnes-type integral
has the following closed form by Lemma~2:
$$
\align
&
\frac1{2\pi i}\int_{-t_0-i\infty}^{-t_0+i\infty}
\frac{\Gamma(-s+t)\,\Gamma(h_0+s+t)\,\Gamma(-t)}
{\Gamma(1+h_0-h_{k+2}+t)}
e^{\pm\pi it}\,\d t
\\ &\qquad
=\frac{\Gamma(-s)}{\Gamma(1-h_{k+2}-s)}
\int_0^1\frac{x^{h_0+s-1}(1-x)^{-h_{k+2}-s}}{(1-x)^{-s}}\,\d x
\allowdisplaybreak &\qquad
=\frac{\Gamma(-s)}{\Gamma(1-h_{k+2}-s)}
\cdot\frac{\Gamma(h_0+s)\,\Gamma(1-h_{k+2})}{\Gamma(1+h_0-h_{k+2}+s)}
\allowdisplaybreak &\qquad
=\frac{\Gamma(h_0+s)\,\Gamma(h_{k+2}+s)\,\Gamma(-s)}
{\Gamma(h_{k+2})\,\Gamma(1+h_0-h_{k+2}+s)}
\cdot\frac{\sin\pi(h_{k+2}+s)}{\sin\pi h_{k+2}}
\\ &\qquad
=\frac{\Gamma(h_0+s)\,\Gamma(h_{k+2}+s)\,\Gamma(-s)}
{\Gamma(h_{k+2})\,\Gamma(1+h_0-h_{k+2}+s)}
\\ &\qquad\quad\times
\biggl(e^{\pi is}\cdot\frac{1-i\cot\pi h_{k+2}}2
+e^{-\pi is}\cdot\frac{1+i\cot\pi h_{k+2}}2\biggr).
\endalign
$$
Substituting this final expression in~\thetag{14} we obtain
$$
\align
&
J_k\biggl(\matrix\format&\,\c\\
h_1, & h_2, & h_3, & \dots, & h_k, & h_{k+1} \\
& 1+h_0-h_3, & 1+h_0-h_4, & \dots, & 1+h_0-h_{k+1}, & 1+h_0-h_{k+2}
\endmatrix\biggr)
\\ &\quad
=\frac{\prod_{j=1}^{k+1}\Gamma(1+h_0-h_j-h_{j+1})}
{\Gamma(h_1)\,\Gamma(h_{k+2})}
\\ &\quad\;\times
\biggl(\frac{1-i\cot\pi h_{k+2}}{4\pi i}
\int_{-s_1-i\infty}^{-s_1+i\infty}
(h_0+2s)\frac{\prod_{j=0}^{k+2}\Gamma(h_j+s)\,\Gamma(-s)}
{\prod_{j=1}^{k+2}\Gamma(1+h_0-h_j+s)}
e^{(\epsilon_{k-1}+1)\pi is}\,\d s
\\ &\quad\;\;
+\frac{1+i\cot\pi h_{k+2}}{4\pi i}
\int_{-s_1-i\infty}^{-s_1+i\infty}
(h_0+2s)\frac{\prod_{j=0}^{k+2}\Gamma(h_j+s)\,\Gamma(-s)}
{\prod_{j=1}^{k+2}\Gamma(1+h_0-h_j+s)}
e^{(\epsilon_{k-1}-1)\pi is}\,\d s
\biggr).
\endalign
$$
If $k$~is even, we take $\epsilon_{k-1}=-1$ in the first integral
and $\epsilon_{k-1}=1$ in the second one. Therefore the both integrals
are equal to
$$
\int_{-s_1-i\infty}^{-s_1+i\infty}
(h_0+2s)\frac{\prod_{j=0}^{k+2}\Gamma(h_j+s)\,\Gamma(-s)}
{\prod_{j=1}^{k+2}\Gamma(1+h_0-h_j+s)}
e^{\epsilon_k\pi is}\,\d s
=2\pi i\cdot F_{k+2}(h_0;h_1,\dots,h_{k+2})
$$
that gives the desired identity~\thetag{4}.
The proof of Theorem is complete.
\enddemo

Another family of multiple integrals
$$
\gather
S(z)
:=\idotsint\limits_{[0,1]^k}
\frac{\prod_{j=1}^kx_j^{a_j-1}(1-x_j)^{b_j-a_j-1}}
{\prod_{i=1}^m(1-zx_1x_2\dotsb x_{r_i})^{c_i}}
\,\d x_1\,\d x_2\dotsb\d x_k,
\tag15
\\
1\le r_1<r_2<\dots<r_m=k,
\endgather
$$
is known due to works of V.~Sorokin~\cite{S1},~\cite{S2}.
Recently, S.~Zlobin~\cite{Zl1},~\cite{Zl2} has proved
(in more general settings)
that the integrals~\thetag{2} can be reduced to the form~\thetag{15}
with $z=1$. Therefore, Theorem gives one a way to reduce
the integrals $S(1)$ to the very-well-poised hypergeometric
series~\thetag{1} under certain conditions on the parameters
$a_j$, $b_j$, $c_i$, and $r_i$ in~\thetag{15}.
In addition, Zlobin~\cite{Zl1} shows that, for integral parameters
in~\thetag{15} satisfying natural restrictions of convergence,
the integral $S(z)$ is a $\Bbb Q[z^{-1}]$-linear combination
of modified multiple polylogarithms
$$
\sum_{n_1\ge n_2\ge\dots\ge n_l\ge1}
\frac{z^{n_1}}{n_1^{s_1}n_2^{s_2}\dotsb n_l^{s_l}}
\qquad\text{with}\quad
s_j\ge1, \; s_j\in\Bbb Z, \; j=1,\dots,l,
$$
where $0\le s_1+s_2+\dots+s_l\le k$ and $0\le l\le m$.

%==================================================

\enddocument
\end